\documentclass[conference]{IEEEtran}
\IEEEoverridecommandlockouts
\usepackage[backend=bibtex,sorting=none]{biblatex}

\usepackage{amsmath,amssymb,amsfonts}
\usepackage{algorithmic}
\usepackage{graphicx}
\usepackage{textcomp}
\usepackage{xcolor}
\usepackage{subcaption}
\usepackage{tikz}
\usepackage{pgfplots}
\usepackage{bm} 
\usepackage{siunitx}
\usepackage[colorlinks = true,
linkcolor = blue,
urlcolor  = blue,
citecolor = blue,
anchorcolor = blue]{hyperref}
\usepackage[capitalise]{cleveref}

\def\BibTeX{{\rm B\kern-.05em{\sc i\kern-.025em b}\kern-.08em
    T\kern-.1667em\lower.7ex\hbox{E}\kern-.125emX}}

\newcommand{\mm}[1]{\textcolor{teal}{#1}}

\begin{document}


\newcommand{\NurbsBasis}{\hat{N}}
\newcommand{\BSplineBasis}{\hat{B}}
\newcommand{\NurbsDegree}{p}
\newcommand{\NumberBasisFunctions}{n}
\newcommand{\NurbsWeight}{w}
\newcommand{\Jac}{\mathbf{J}_F}
\newcommand{\GeometryBasis}{G}

\newcommand\st{\mathrm{st}}
\newcommand\rt{\mathrm{rt}}
\newcommand\ag{\mathrm{ag}}
\newcommand\source{\mathrm{src}}

\newcommand{\der}{\operatorname{d}}
\newcommand{\MagVecPot}{\mathbf{A}}
\newcommand{\MagVecPotz}{A_z}
\newcommand{\MagVecPotzRt}{A_{z,\rt}}
\newcommand{\MagVecPotzSt}{A_{z,\st}}
\newcommand{\CurrentDensity}{\mathbf{J}}
\newcommand{\CurrentDensityTwoD}{J}

\newcommand{\MagFluxDensity}{\mathbf{B}}
\newcommand{\Remanence}{\mathbf{B_\mathrm{r}}}
\newcommand{\RemanenceRed}{\mathbf{B^\bot_\mathrm{r}}}
\newcommand{\RotAngle}{\beta}
\newcommand{\MagnetAngle}{\alpha}
\newcommand{\CurrentAngle}{\varphi_0}

\newcommand{\ScaleR}{k_\mathrm{R}}
\newcommand{\LossJoule}{P_\mathrm{J}}


\newcommand{\ansatzFunction}{N} 
\newcommand{\ansatzFunctionIndex}{j}
\newcommand{\testSymbol}{v}
\newcommand{\testFunction}{N} 
\newcommand{\testFunctionIndex}{i}

\newcommand{\couplingSymbol}{G}
\newcommand{\couplingMatrix}{\mathbf{\couplingSymbol}}
\newcommand{\stiffnessSymbol}{K}
\newcommand{\stiffnessMatrix}{\mathbf{\stiffnessSymbol}}

\newcommand{\solutionSymbol}{u}
\newcommand{\solutionVector}{\mathbf{\solutionSymbol}}
\newcommand{\rhsSymbol}{b}
\newcommand{\rhsVector}{\mathbf{\rhsSymbol}}
\newcommand{\optiSymbol}{x}
\newcommand{\optiVector}{\mathbf{\optiSymbol}}
\newcommand{\stateSymbol}{e}
\newcommand{\stateVector}{\mathbf{\stateSymbol}}

\newcommand{\mortarSymbol}{\lambda}
\newcommand{\mortarVector}{\bm{\mortarSymbol}}

\newcommand{\opt}{\mathrm{opt}}
\newcommand{\fopt}{f_\opt}
\newcommand{\RotMat}{\mathbf{R}_\RotAngle}
\newcommand{\ControlPointSymbol}{C}
\newcommand{\OptiCtrlPoint}{\mathbf{\ControlPointSymbol}}
\newcommand{\ParameterSymbol}{P}
\newcommand{\OptiParameter}{\mathbf{\ParameterSymbol}}
\newcommand{\adjointSymbol}{\gamma}
\newcommand{\adjointVector}{\bm{\adjointSymbol}}
\newcommand{\Torque}{T}
\newcommand{\MeanTorque}{\overline{\Torque}}
\newcommand{\StdTorque}{\hat{\Torque}}
\newcommand{\Amagnet}{A_\mathrm{Magnet}}
\newcommand{\Ttarget}{\Torque_\mathrm{Target}}

\newcommand{\OmegaRotor}{\Omega_\rt}
\newcommand{\OmegaStator}{\Omega_\st}
\newcommand{\GammaAirGap}{\Gamma_\ag}
\newcommand\IntegG{\mathrm{d}\Gamma}
\newcommand\Integ{\mathrm{d}\Omega}

\newcommand{\Iapp}{I_{\mathrm{app}}}
\newcommand{\nWind}{n_{\mathrm{wind}}}
\newcommand{\Acoil}{A_{\mathrm{coil}}}
\newcommand{\polePair}{p}

\definecolor{TUDa-2a}{HTML}{009CDA}
\definecolor{TUDa-2b}{HTML}{0083CC}
\definecolor{TUDa-3d}{HTML}{0071F3}
\definecolor{TUDa-3a}{HTML}{50B695}
\definecolor{TUDa-9b}{HTML}{E6001A}

\title{Spline-Based Rotor and Stator Optimization of a Permanent Magnet Synchronous Motor}

\author{\IEEEauthorblockN{1\textsuperscript{st} Michael Wiesheu}
\IEEEauthorblockA{\textit{Computational Electromagnetics Group} \\
\textit{Technical University of Darmstadt}\\
Darmstadt, Germany \\
michael.wiesheu@tu-darmstadt.de}
\and
\IEEEauthorblockN{2\textsuperscript{nd} Theodor Komann}
\IEEEauthorblockA{\textit{Department of Mathematics} \\
\textit{Technical University of Darmstadt}\\
Darmstadt, Germany  \\
komann@mathematik.tu-darmstadt.de}
\and
\IEEEauthorblockN{3\textsuperscript{rd} Melina Merkel}
\IEEEauthorblockA{\textit{Computational Electromagnetics Group} \\
\textit{Technical University of Darmstadt}\\
Darmstadt, Germany  \\
melina.merkel@tu-darmstadt.de}
\and
\IEEEauthorblockN{4\textsuperscript{th} Sebastian Schöps}
\IEEEauthorblockA{\textit{Computational Electromagnetics Group} \\
\textit{Technical University of Darmstadt}\\
Darmstadt, Germany  \\
sebastian.schoeps@tu-darmstadt.de}
\and
\IEEEauthorblockN{5\textsuperscript{th} Stefan Ulbrich}
\IEEEauthorblockA{\textit{Department of Mathematics} \\
\textit{Technical University of Darmstadt}\\
Darmstadt, Germany  \\
ulbrich@mathematik.tu-darmstadt.de}
\and
\IEEEauthorblockN{6\textsuperscript{th} Idoia Cortes Garcia}
\IEEEauthorblockA{\textit{Department of Mechanical Engineering} \\
\textit{Eindhoven University of Technology}\\
Eindhoven, The Netherlands \\
i.cortes.garcia@tue.nl}
}

\maketitle

\begin{abstract}

This work features the optimization of a Permanent Magnet Synchronous Motor using 2D nonlinear simulations in an Isogeometric Analysis framework. The rotor and stator designs are optimized for both geometric parameters and surface shapes via modifications of control points. The scaling laws for magnetism are employed to allow for axial and radial scaling, enabling a thorough optimization of all critical machine parameters for multiple operating points. The process is carried out in a gradient-based fashion with the objectives of lowering motor material cost, torque ripple and losses. It is shown that the optimization can be efficiently conducted for many optimization variables and all objective values can be reduced.
\end{abstract}


\section{Introduction}

Optimizing electric machines is a challenging task. On the one hand, there are a variety of options to change the geometric design, e.g., parameters (such as magnet width, yoke thickness or motor length) or shape adjustments (such as the rotor and stator surfaces). This can lead to a high dimensional design space. On the other hand, there are multiple quantities of interest which should be improved, such as the motor cost, performance (given quantified e.g. by the motor torque), or efficiency.

This has led to a variety of published research articles and surveys that feature Finite Element (FE) based optimization of e.g. cost \cite{Meddour_2023xx} or torque ripple \cite{Sanchez_2022aa} with methods such as parameter optimization \cite{Omar_2022xx}, shape optimization \cite{Merkel_2021ab} or topology optimization \cite{Nishanth_2022xx}. Gradient-free algorithms, like particle swarms or genetic algorithms, are commonly applied in these studies due to their simplicity and promise of achieving a global optimum. Yet, these approaches can become computationally expensive when the optimization has a large design space, i.e., there are many variables to optimize. 

This work extends the previously proposed method in \cite{Wiesheu_2023ab}, where the rotor of a Permanent Magnet Synchronous Machine (PMSM) is optimized for both parameters and shape with gradient based optimization. Here,  both, the rotor and stator parameters and surfaces are optimized, and multiple operating points are considered for a different motor geometry. A snapshot of the code can be found in \cite{michaelwiesheu_2024_10726571}.

\section{Methodology}
This section briefly describes the underlying numerical methods and optimization techniques. It focuses on new features that are not covered in the previous work \cite{Wiesheu_2023ab}, which contains methodological explanations in detail. For more information, the interested reader is also referred to preliminary work deriving the mathematical foundations \cite{Egger_2022ab,Bontinck_2018ac}.

\subsection{Numerical modeling}
To describe the electromagnetic behavior of the electric machine, we employ the magnetostatic problem
\begin{equation}
    \nabla\times\left(\nu\nabla\times\MagVecPot\right)= \CurrentDensity+\nabla\times\left(\nu\Remanence\right),
    \label{eq:Maxwell1}
\end{equation}
with the magnetic vector potential $\MagVecPot$, the (nonlinear) reluctivity $\nu$, the current source density $\CurrentDensity$ and the magnet remanence $\Remanence$. The magnetic flux density $\MagFluxDensity$ is defined via $\MagFluxDensity = \nabla\times\MagVecPot$. This formulation is common for PMSM machines, where eddy currents are neglected due to the lamination of the iron cores \cite{Salon_1995aa}. Reducing \eqref{eq:Maxwell1} to 2D and considering rotor and stator domains $\OmegaRotor$ and $\OmegaStator$ separately yields the Poisson problem
\begin{equation}
    \begin{cases}
    \nabla \cdotp ( \nu \nabla \MagVecPotzRt) =\nu \nabla \cdotp \RemanenceRed & \mathrm{in}\ \OmegaRotor\\
    \nabla \cdotp ( \nu \nabla \MagVecPotzSt) =-\CurrentDensityTwoD & \mathrm{in}\ \OmegaStator,
    \end{cases} \label{eq:StrongRotorStator}
\end{equation}
where $\MagVecPot = \begin{pmatrix}
0 & 0 & \MagVecPotz
\end{pmatrix}^\top$. Homogeneous Dirichlet boundaries are applied at $\Gamma_\mathrm{d}$, antiperiodic ones at $\Gamma_\mathrm{ap}$ and the coupling of rotor and stator is performed at the air gap interface $\Gamma_\mathrm{ag}$, see \cref{fig:PMSMgeometry}. 3D effects, e.g. due to the end windings, will be neglected.
\begin{figure*}[h]
    \centering
    \input{PMSMgeometry.tikz}
    \caption{Parametrization of the PMSM geometry including parameter names, material definitions and boundary conditions. The model is originally based on \cite{Pahner_1998ab}.}
    \label{fig:PMSMgeometry}
\end{figure*}
The right-hand side in \eqref{eq:StrongRotorStator} is given by the 2D components of the remanence and the three phase current. The 2D remanence is given by $\RemanenceRed = B_{\mathrm{r}}\begin{pmatrix}
-\sin( \MagnetAngle ) & \cos( \MagnetAngle )
\end{pmatrix}^\top$ defined by the remanence $B_{\mathrm{r}}$ and orientation $\MagnetAngle$ and the three phase current is given by 
$\CurrentDensityTwoD=\sum\nolimits_{k}\CurrentDensityTwoD^{(k)}$ 
with
\begin{align}
\CurrentDensityTwoD^{(k)} & = J_0\sin\left(\polePair\RotAngle +\CurrentAngle +\frac{2\pi }{3} k\right),\label{eq:CurrentDefinition}
\end{align}
where $k\in \{0,1,2\}$ indicates the $k$-th phase. 
For synchronous operation, this can be expressed by a current density $J_0$, the pole pair number $p$, the electric phase offset $\CurrentAngle$ and the rotation angle $\RotAngle$.

The PMSM considered in this work is illustrated in \cref{fig:PMSMgeometry}. Due to symmetry, only one sixth of the motor is simulated. The rotor and stator iron cores are shown in gray, the (homogenized) copper slots in red, the rotor magnet in green and the air gap and air pockets in blue. In addition, the parameter names for later optimization are given. 

Discretization of \eqref{eq:StrongRotorStator} with Isogeometric Analysis (IGA) and coupling of the rotor and stator surface with Harmonic Mortaring \cite{Bontinck_2018ac} yields the matrix system
\begin{equation}
    \underbrace{\begin{pmatrix}
    \stiffnessMatrix_\rt & \mathbf{0} & -\couplingMatrix_\rt \\
    \mathbf{0} & \stiffnessMatrix_\st & \couplingMatrix_\st\RotMat \\
    -\couplingMatrix_\rt^\top & \RotMat^\top\couplingMatrix_\st^\top & \mathbf{0}
    \end{pmatrix}}_{=:\stiffnessMatrix}\underbrace{\begin{pmatrix}
    \solutionVector_\rt \\
    \solutionVector_\st \\
    \mortarVector
    \end{pmatrix}}_{=:\solutionVector} =\underbrace{\begin{pmatrix}
    \rhsVector_\rt \\
    \rhsVector_\st \\
    \mathbf{0}
    \end{pmatrix}}_{=:\rhsVector} \label{eq:CoupledMatrixSystem}
\end{equation}
with the stiffness matrices $\stiffnessMatrix_\rt$ and $\stiffnessMatrix_\st$, coupling matrices $\couplingMatrix_\rt$ and $\couplingMatrix_\st$, the rotation matrix $\RotMat$ depending on the rotation angle $\RotAngle$, the magnetic potential $\solutionVector_\rt$ and $\solutionVector_\st$, the Lagrange multipliers $\mortarVector$ and the right-hand side $\rhsVector_\rt$ and $\rhsVector_\st$. Because of the Harmonic Mortaring, only the rotation matrix $\RotMat$, which contains sine and cosine function values, needs to be recalculated for different $\RotAngle$. For conciseness, \eqref{eq:CoupledMatrixSystem} is rewritten as the state equation
\begin{equation}
\stateVector(\solutionVector) = \stiffnessMatrix(\solutionVector)\solutionVector - \rhsVector, \label{eq:CoupledMatrixSystemred}
\end{equation}
where $\stiffnessMatrix$ depends also on $\solutionVector$ due to material nonlinearities. 
After solving \eqref{eq:CoupledMatrixSystemred} with a Newton scheme, the electromagnetic torque $T_\RotAngle$ of the motor is determined with
\begin{equation}
    T_\RotAngle(\solutionVector) = -\solutionVector_\st^\top \couplingMatrix_\st \RotMat^\prime \mortarVector L \ScaleR^2 \label{eq:Torque}
\end{equation}
by introducing the motor length $L$ and the radial scaling $\ScaleR$.

Evaluating \eqref{eq:Torque} for multiple $\RotAngle \in \{\RotAngle_1,\RotAngle_2,...,\RotAngle_N \}$, we can calculate the average torque $\MeanTorque$ and the torque standard deviation $\StdTorque$ with
\begin{align}
    \MeanTorque &=\frac{1}{N}\sum\nolimits _{\RotAngle } T_{\RotAngle } & \StdTorque &=\sqrt{\frac{1}{N}\sum\nolimits _{\RotAngle }( T_{\RotAngle } -\MeanTorque)^{2}}.
    \label{eq:Optimization:Torques}
\end{align}
The average torque is linked to the motor power and the torque standard deviation influences the rotation smoothness, so both values are relevant quantities to be optimized.

\subsection{Scaling laws for the Finite Element solution}
For convenience in the following optimization steps, scaling laws of the FE solution are introduced \cite{Nell_2019aa,Stipetic_2015xx}. This allows for the inclusion of the motor's radial scaling (given by $k_\mathrm{R}$) and axial length (given by $L$) as optimization variables. \cref{tab:ScalingLaws} summarizes the necessary variables and their respective scaling with respect to $k_\mathrm{R}$ and $L$.
\begin{table}[!h]
    \centering
   \renewcommand{\arraystretch}{1.2}
    \caption{Scaling laws with respect to $\ScaleR$ and $L$, according to \cite{Nell_2019aa}.}
    \begin{tabular}{|l|c|c|}
    \hline 
    \textbf{Parameter} & \textbf{Variable} & $\displaystyle \varpropto $ \\
    \hline 
    Magnetic flux density & $\displaystyle \textbf{B}$ & $\displaystyle 1$ \\
    Magnetic field strength & $\displaystyle \textbf{H}$ & $\displaystyle 1$ \\
    Cross sectional area & $A$ & $\displaystyle k_{\mathrm{R}}^{2}$ \\
    Stator current & $\displaystyle I$ & $\displaystyle k_{\mathrm{R}}$ \\
    Stator current density & $\displaystyle J$ & $\displaystyle k_{\mathrm{R}}^{-1}$ \\
    Torque & $\displaystyle T$ & $\displaystyle k_{\mathrm{R}}^{2} L$ \\
    Material cost & $\displaystyle M$ & $\displaystyle k_{\mathrm{R}}^{2} L$ \\
    \hline
    \end{tabular}
    \label{tab:ScalingLaws}
\end{table}

Note that the stator current $I$ must increase by $\ScaleR$ in order for the magnetic field to remain unchanged if radial scaling $\ScaleR$ is applied. That implies a scaling of the current density $J$ by $\ScaleR^{-1}$, which can lead to unrealistic results, as the maximum current density can not be increased arbitrarily. Therefore, we scale the current density from \eqref{eq:CurrentDefinition} by $\ScaleR$, such that $J_0$ in the final configuration always remains the same. 

For the optimization, we further calculate the motor material cost $M$ by
\begin{equation}
M =Lk_{R}^{2}\sum\nolimits_i \rho _{i} A_{i} c_{i}
\end{equation}
with the material density $\rho_i$, cross sectional $A_i$ and cost per mass $c_i$ for the i-th material, respectively. In addition to $M$ we define the Joule stator losses $\LossJoule$ \cite{Bianchi_2005aa} for one slot as
\begin{equation}
    \LossJoule =CA_{\mathrm{slot}} L,
    \label{eq:LossJoule}
\end{equation}
where $C$ is a constant taking into account the stator current density, conductivity and number of windings. This is the only loss mechanism considered here and is used to penalize less efficient designs.

For the necessary derivatives in the gradient based optimization, the simple relationships 
\begin{align}
\frac{\mathrm{d} T}{\mathrm{d} L} & =\frac{T}{L} & \frac{\mathrm{d} T}{\mathrm{d} k_{\mathrm{R}}} & =2\frac{T}{k_{\mathrm{R}}}  \nonumber \\ 
\frac{\mathrm{d} M}{\mathrm{d} L} & =\frac{M}{L} & \frac{\mathrm{d} M}{\mathrm{d} k_{\mathrm{R}}} & =2\frac{M}{k_{\mathrm{R}}} \label{eq:DerScale} \\
\frac{\mathrm{d} P_{\mathrm{J}}}{\mathrm{d} L} & =\frac{P_{\mathrm{J}}}{L} & \frac{\mathrm{d} P_{\mathrm{J}}}{\mathrm{d} k_{\mathrm{R}}} & =2\frac{P_{\mathrm{J}}}{k_{\mathrm{R}}} \nonumber
\end{align}
can be exploited in addition to the derivatives given in \cite{Wiesheu_2023ab} to include the radial and axial dimensions. 

\subsection{Optimization}
The objective function consists of three components that should be minimized: the motor mass $M$, the torque ripple $\StdTorque$ and the Joule losses $\LossJoule$. We employ a multi-objective optimization approach with a weighted sum, i.e., every objective is scaled by a weight $m_i$. The torque ripple is evaluated for different operating points and summed up in  $\StdTorque$. Here, we choose the operating points in the set $\CurrentDensityTwoD \in J_0\cdotp\{0, 0.5, 1\}$ and $\RotAngle \in \{0,2,...,18\}$ reducing both cogging torque and the torque ripple.

As constraints, first, the magnetostatic formulation \eqref{eq:CoupledMatrixSystemred} must be fulfilled. Second, the mean torque $\MeanTorque$ must meet the target value $T_{\mathrm{Target}}$. Third, geometric constraints $\mathbf{g}$ are applied, such that there are no overlaps between different domains. Overall, the optimization problem is written as
\begin{equation}
    \begin{array}{l}
        \min \fopt(\optiVector,\solutionVector) = m_1 M + m_2 \StdTorque + m_3 \LossJoule \\[2.5pt]
        \mathrm{s.t.} \\ 
        \begin{cases}
        \hspace{0.5cm} \stateVector(\optiVector,\solutionVector) = \mathbf{0} & \mathrm{State~equation}\\
        \hspace{0.5cm} \MeanTorque \geq T_{\mathrm{Target}} & \mathrm{Fulfill~target~torque} \\
        \hspace{0.5cm} \mathbf{g}(\optiVector) \leq \mathbf{0} & \mathrm{Geometric~feasibility}. 
        \end{cases}
    \end{array}\label{eq:OptimizationProblem}
\end{equation}

The design vector $\optiVector$ contains the parameters shown in \cref{fig:PMSMgeometry}, the variables $\CurrentAngle$, $L$, $\ScaleR$ as well as the radial offsets of the rotor and stator surface control points. 

The derivatives of \eqref{eq:OptimizationProblem} are calculated with \eqref{eq:DerScale} as well as the analytical derivatives from \cite{Wiesheu_2023ab} using the adjoint method. This allows for an efficient optimization process despite the large number of design variables.

\section{Results}
The optimization is carried out with \textit{MATLAB}\textsuperscript{\textregistered}'s solver \textit{fmincon} using the \textit{interior-point} method in an IGA framework with \textit{GeoPDEs} \cite{Vazquez_2016aa}. The material cost is chosen as a dimensionless number with $c_\mathrm{Iron} = 2$, $c_\mathrm{Copper} = 10$ and $c_\mathrm{Magnet} = 50$ to roughly represent the proportion of real material cost. The objective weights in \cref{eq:OptimizationProblem} are set to $m_1 = 0.05$ $m_2 = 10$ $m_3 = 4$ such that the different objectives are of comparable magnitude. The maximum current density is set to $\displaystyle J_0=\SI{3.2}{\ampere\per\mm\squared}$. As torque requirement, a torque of $T_\mathrm{Target} \geq \SI{1.5}{Nm}$ must be maintained. For the iron, the nonlinear BH-characteristics of M330-50A (M27) is used \cite{FEMM}. \cref{tab:Parameters} summarizes the initial parameter values with lower and upper bounds.
\begin{table}[h]
    \centering
    \caption{Motor parameters and bounds for the optimization.}
    \begin{tabular}{|l|r|r|r|r|}
    \hline
    \textbf{Parameter name} & \textbf{initial} & \textbf{min} & \textbf{max} & \textbf{opt} \\
    \hline
    $L$ & \SI{100}{\mm} & \SI{80}{\mm} & \SI{120}{\mm} & \SI{80}{\mm} \\
    $\CurrentAngle$ & \SI{0}{\deg} & \SI{-20}{\deg} & \SI{20}{\deg} & \SI{-11.6}{\deg} \\
    $\ScaleR$ & 1 & 0.5 & 2 & 1.048 \\
    MAG & \SI{7}{\mm} & \SI{6}{\mm} & \SI{15}{\mm} & \SI{6.10}{\mm} \\
    MH & \SI{7}{\mm} & \SI{2}{\mm} & \SI{12}{\mm} & \SI{2.74}{\mm} \\
    MW & \SI{19}{\mm} & \SI{10}{\mm} & \SI{25}{\mm} & \SI{20.73}{\mm} \\
    SD1 & \SI{135}{\mm} & \SI{90}{\mm} & \SI{160}{\mm} & \SI{130.29}{\mm} \\
    SR1 & \SI{1}{\mm} & \SI{0.5}{\mm} & \SI{2}{\mm} & \SI{2.0}{\mm} \\
    SW1 & \SI{4}{\mm} & \SI{2}{\mm} & \SI{6}{\mm} & \SI{4.06}{\mm} \\
    SW2 & \SI{2.3}{\mm} & \SI{1}{\mm} & \SI{4}{\mm} & \SI{1.0}{\mm}\\
    SW3 & \SI{1}{\mm} & \SI{0.5}{\mm} & \SI{1.5}{\mm} & \SI{0.50}{\mm} \\
    SW4 & \SI{8.25}{\mm} & \SI{2}{\mm} & \SI{20}{\mm} & \SI{5.43}{\mm} \\
    \hline
    \end{tabular}
    \label{tab:Parameters}
\end{table}

\begin{figure*}[t]
    \centering
    \begin{subfigure}[t]{.49\linewidth}
      \centering
    \includegraphics[width=\linewidth]{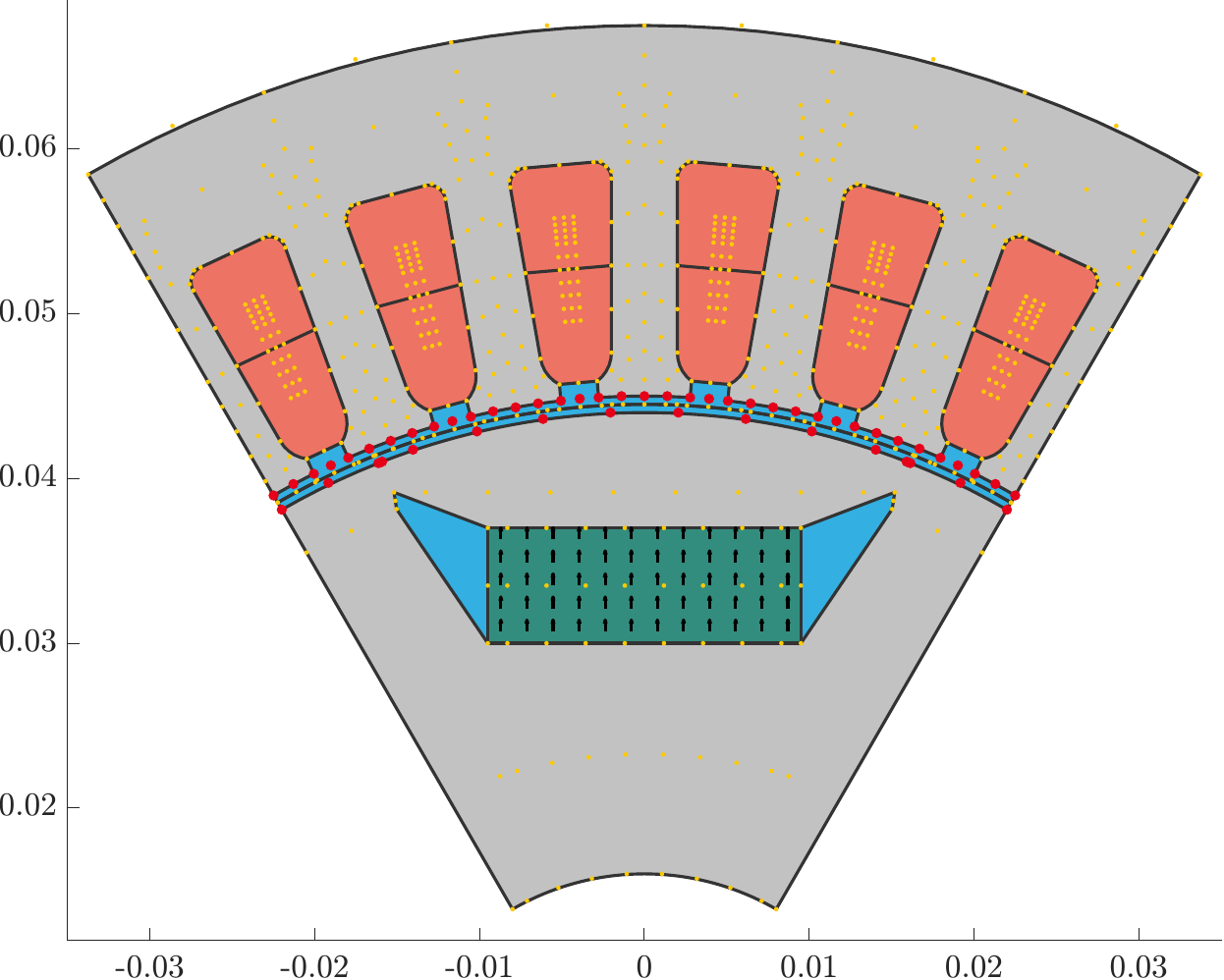}
        \caption{Geometry representation of the initial motor geometry.}
      \label{fig:GeometryInit}
    \end{subfigure} \hfill
    \begin{subfigure}[t]{.49\linewidth}
      \centering
        \includegraphics[width=\linewidth]{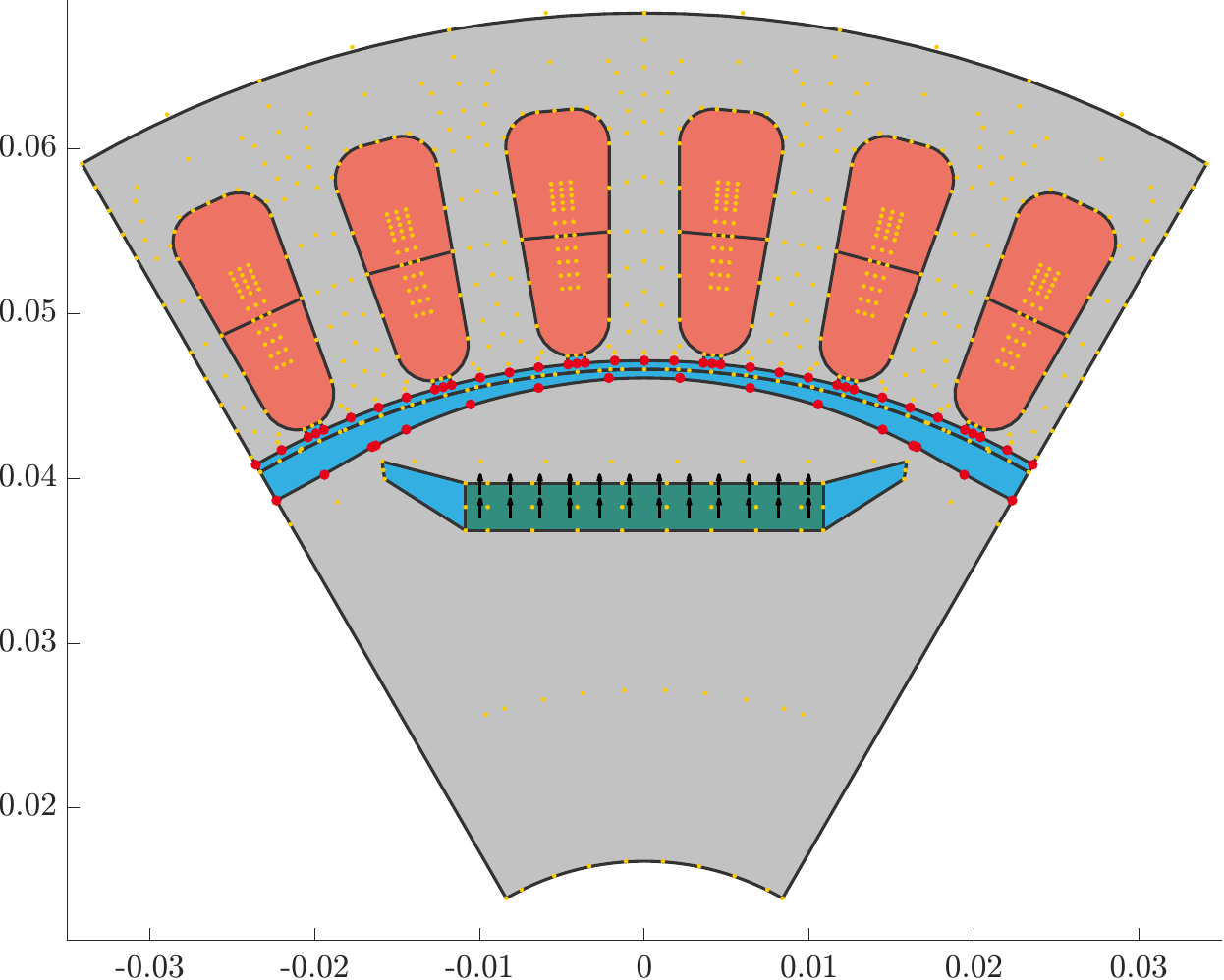}
      \caption{Geometry representation of the optimized motor geometry.}
      \label{fig:GeometryOpt}
    \end{subfigure}
    \caption{Comparison of initial and optimized geometry. The yellow points show the control points of the motor, the red control points may move during the optimization.}
    \label{fig:ComparisonGeometry}
\end{figure*}

\begin{figure*}[t]
    \centering
    \begin{subfigure}[t]{.49\linewidth}
      \centering
    \includegraphics[width=\linewidth]{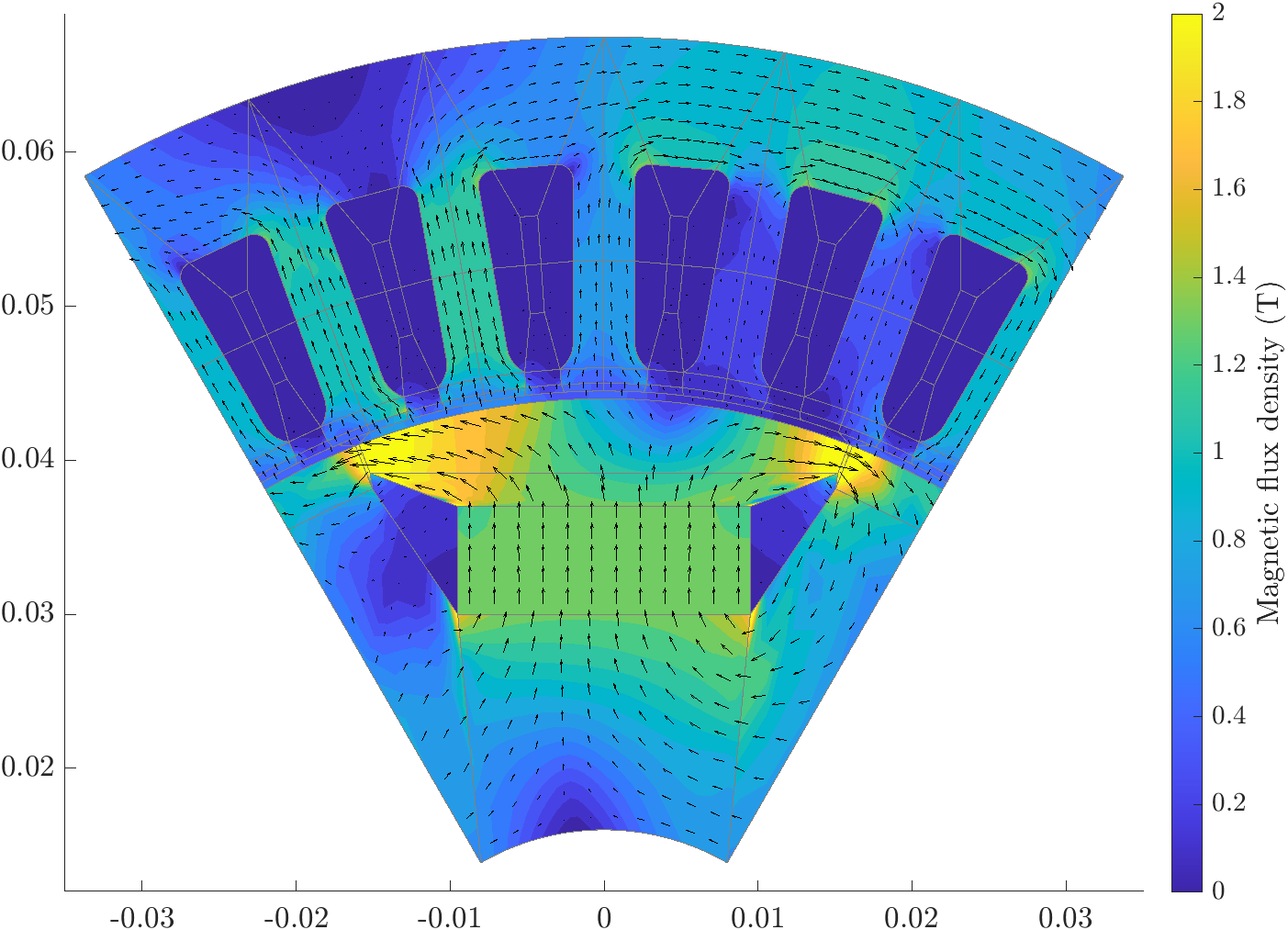}
        \caption{Magnetic flux density of the original motor.}
      \label{fig:FieldInit}
    \end{subfigure} \hfill
    \begin{subfigure}[t]{.49\linewidth}
      \centering
        \includegraphics[width=\linewidth]{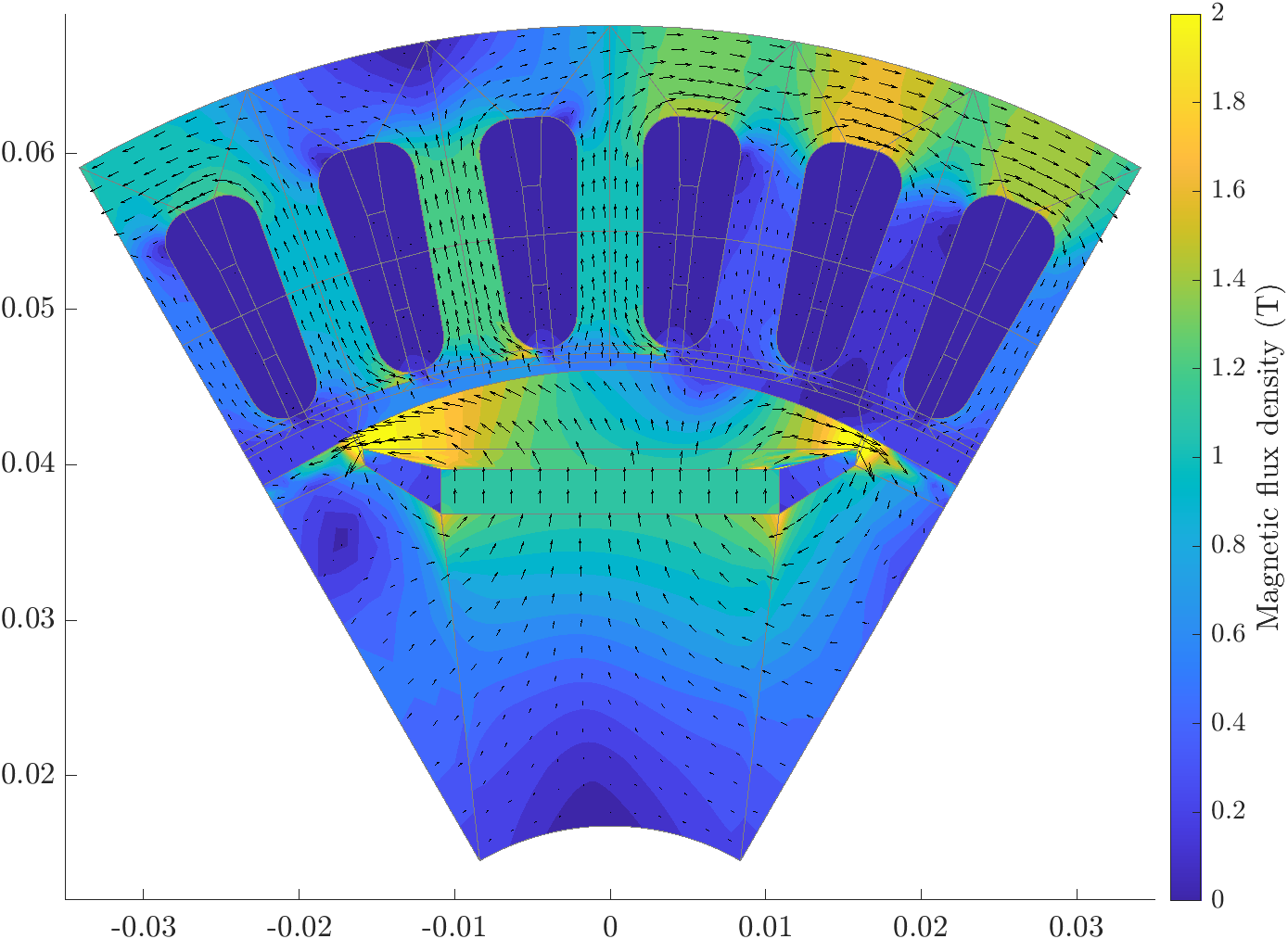}
      \caption{Magnetic flux density of the optimized motor.}
      \label{fig:FieldOpt}
    \end{subfigure}
    \caption{Comparison of the magnetic flux density. }
    \label{fig:ComparisonField}
\end{figure*}

Overall, there are 22 design variables, where 12 come from the parameters and 10 from the (symmetrically chosen) control points. The FE  simulation is carried out with 4095 degrees of freedom and 741 control points. 
These are highlighted in \cref{fig:GeometryInit}, where the initial and geometry is shown. The magnetic flux density for the initial design at $\RotAngle=0^\circ$ is shown in \cref{fig:FieldInit}.


The optimization is carried out on a 8-core laptop\footnote{Intel\textsuperscript{\textregistered} Core\textsuperscript{\texttrademark} i7-1165G7@2.8GHz with \SI{16}{GB} RAM.} in 120 iterations within six hours. The results of the optimization are presented in \cref{tab:Results}. The optimized geometry is given in \cref{fig:GeometryOpt}. Two features of the new geometry are that the magnet size is greatly reduced, leading to a lower cost and that the rotor surface is not circular any more, which reduces the torque ripple. As there is a trade off between magnet size (influences cost) and the coil size (influences loss), the weights in \eqref{eq:OptimizationProblem} can be further adjusted to generate multiple designs in a Pareto front. The distribution of the magnetic flux density for the optimized design is shown in \cref{fig:FieldOpt}.  

\begin{table}[!h]
        \centering
    \caption{Comparison of objectives for the initial and optimized motor.}
    \label{tab:Results}
    \begin{tabular}{c|c|c|c|}
    \cline{2-4} 
      & \textbf{Initial} & \textbf{Optimized} & \textbf{Change} \\
    \hline 
     \multicolumn{1}{|l|}{$f_\mathrm{opt}$} & 2.0555 & 0.5743 & -72.1\% \\
    \hline 
     \multicolumn{1}{|l|}{Cost Motor} & 11.1506 & 7.2520 & -35.0\% \\
    \hline 
     \multicolumn{1}{|l|}{ \ \ \ Cost Iron} & 2.534 & 2.1006 & -17.1\% \\
    \hline 
     \multicolumn{1}{|l|}{ \ \ \ Cost Copper} & 3.6291 & 3.2774 & -9.69\% \\
    \hline 
     \multicolumn{1}{|l|}{ \ \ \ Cost Magnet} & 4.9875 & 1.8739 & -62.4\% \\
    \hline 
     \multicolumn{1}{|l|}{Torque Ripple at OP} & 0.129 & 0.0013 & -99.9\% \\
    \hline 
     \multicolumn{1}{|l|}{Power Loss} & 0.0520 & 0.0497 & -4.42\% \\
     \hline
    \end{tabular}
\end{table}

\begin{figure*}[t]
    \centering
    \begin{subfigure}[t]{.49\linewidth}
        \centering
%
%
\definecolor{mycolor1}{rgb}{0.36471,0.52157,0.76471}%
\definecolor{mycolor2}{rgb}{0.31373,0.71373,0.58431}%
\begin{tikzpicture}

\begin{axis}[%
width=\linewidth,
height=5.5cm,
xmin=0,
xmax=20,
xlabel style={font=\color{white!15!black}},
xlabel={$\text{Rotation angle (}^\circ\text{)}$},
ymin=-0.2,
ymax=1.8,
ylabel style={yshift=-0.3cm,font=\color{white!15!black}},
ylabel={Torque (Nm)},
xtick={0,2,...,20},
legend pos=south east,
legend cell align={left},
axis background/.style={fill=white},
axis x line*=bottom,
axis y line*=left,
xmajorgrids,
ymajorgrids
] 
\addplot [color=mycolor1, line width=1.5pt, ,legend image post style={mark={*}}]
table[x index=0,y index= 1,col sep=comma] {TorqueData.csv};

\addplot [color=mycolor1, line width=1.5pt, forget plot]
table[x index=0,y index= 2,col sep=comma] {TorqueData.csv}; 

\addplot [color=mycolor1, line width=1.5pt, forget plot]
table[x index=0,y index= 3,col sep=comma] {TorqueData.csv}; 

\addplot [only marks, mark=*, mark options={}, mark size=2pt, color=mycolor1, fill=mycolor1, forget plot]
table[x index=0,y index= 1,col sep=comma] {TorqueScatter.csv}; 

\addplot [only marks, mark=*, mark options={}, mark size=2pt, color=mycolor1, fill=mycolor1, forget plot]
table[x index=0,y index= 2,col sep=comma] {TorqueScatter.csv};

\addplot [only marks, mark=*, mark options={}, mark size=2pt, color=mycolor1, fill=mycolor1, forget plot]
table[x index=0,y index= 3,col sep=comma] {TorqueScatter.csv};
\addlegendentry{Initial}

\addplot [color=mycolor2, line width=1.5pt, legend image post style={mark={square*}}]
table[x index=0,y index= 4,col sep=comma] {TorqueData.csv};

\addplot [color=mycolor2, line width=1.5pt, forget plot]
table[x index=0,y index= 5,col sep=comma] {TorqueData.csv}; 

\addplot [color=mycolor2, line width=1.5pt, forget plot]
table[x index=0,y index= 6,col sep=comma] {TorqueData.csv}; 

\addplot [only marks, mark=square*, mark options={}, mark size=2pt, color=mycolor2, fill=mycolor2 , forget plot]
table[x index=0,y index= 4,col sep=comma] {TorqueScatter.csv}; 

\addplot [only marks, mark=square*, mark options={}, mark size=2pt, color=mycolor2, fill=mycolor2 , forget plot]
table[x index=0,y index= 5,col sep=comma] {TorqueScatter.csv};

\addplot [only marks, mark=square*, mark options={}, mark size=2pt, color=mycolor2, fill=mycolor2 , forget plot]
table[x index=0,y index= 6,col sep=comma] {TorqueScatter.csv};
\addlegendentry{Optimized}

\end{axis}
\end{tikzpicture}%
        \caption{Torque profiles of the initial and optimized motor for $\CurrentDensityTwoD \in J_0\cdotp\{0, 0.5, 1\}$. The markers highlight the evaluated angles in the optimization process.}
        \label{fig:TorqueProfiles}
    \end{subfigure} \hfill
    \begin{subfigure}[t]{.49\linewidth}
            \centering
        \includegraphics[width=\linewidth]{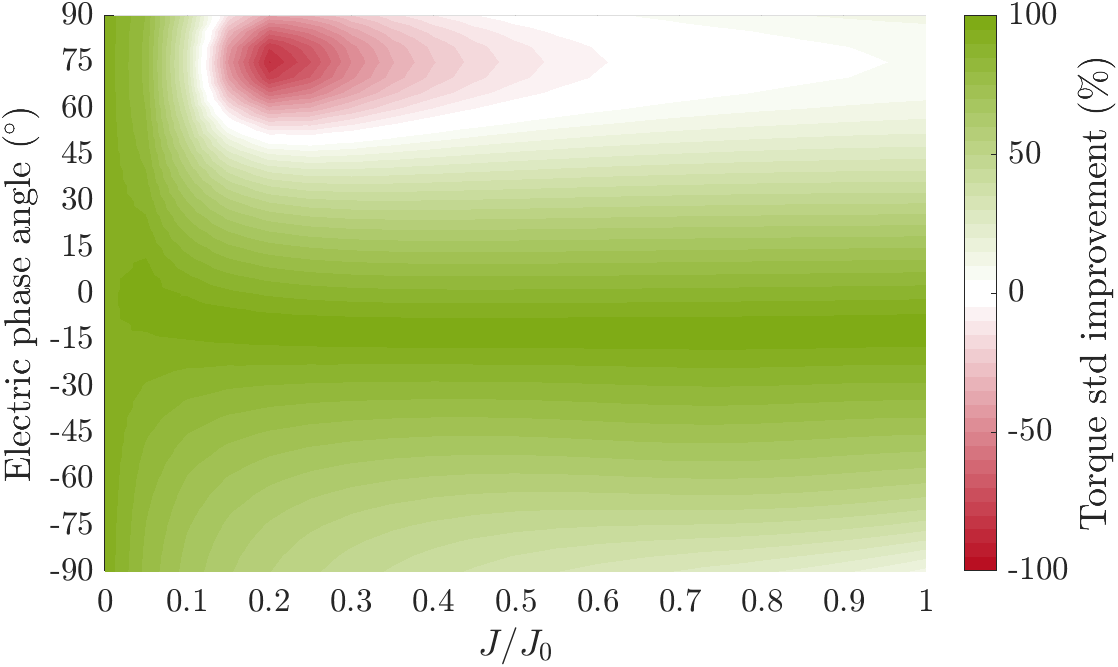}
        \caption{Percentage change of torque ripple plotted over the electric phase offset and the phase current. The torque ripple decreases significantly over a large range of operating points.}
        \label{fig:TorqueMap}
    \end{subfigure}
    \caption{Comparison of initial and optimized torque profiles for different operating points.}
    \label{fig:TorqueResults}
\end{figure*}
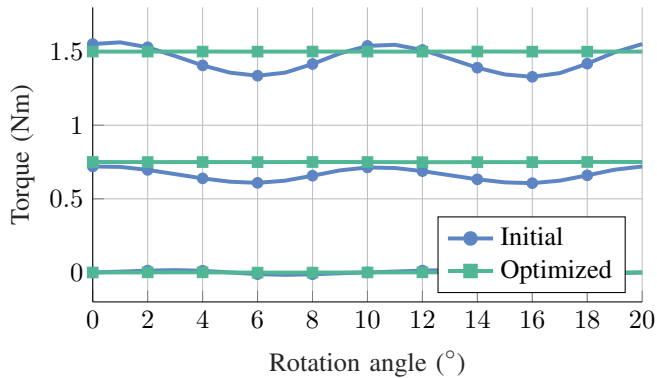
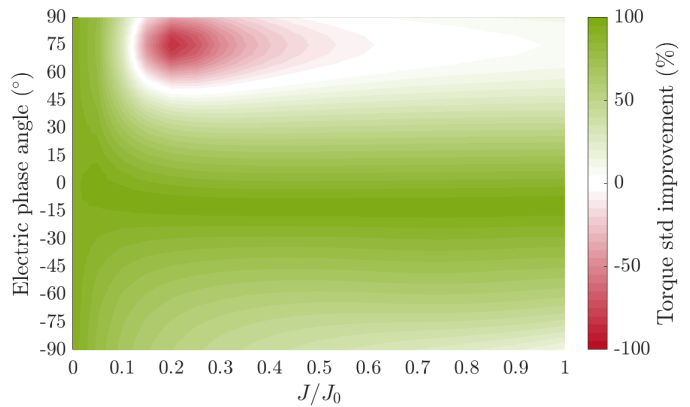

The findings can be summarized as follows:
\begin{itemize}
    \item The motor cost, torque ripple and power loss can be reduced. The decrease is \SI{35.0}{\percent}, \SI{99.0}{\percent} and \SI{4.42}{\percent} respectively.
    \item Compared to the initial design (\cref{fig:FieldInit}), the magnetic flux density in the optimized design (\cref{fig:FieldOpt}) is increased in the stator teeth and yoke. The optimizer chooses the teeth and yoke width such that the saturation of the iron is optimally exploited.
    \item In contrast to the rotor surface, the stator surface remains almost unchanged and circular. This means that adapting the rotor surface has the most relevant influence to reduce the torque ripple.
    \item The thickness of the iron bridges in the rotor are reduced. This allows for a smaller magnet size but increases mechanical stresses at the bridges.
    \item The optimization sets the motor length to the minimum value and compensates this with increased radial scale. In order to find the best motor length, mechanical stress constraints need to be added in the future such that the radial scale is not increased arbitrarily.
\end{itemize}

A comparison of the torque profiles for the initial and optimized motor is found in \cref{fig:TorqueResults}. In \cref{fig:TorqueProfiles} the initial and optimized torque profiles for the three operating points that were used for optimization are compared. \cref{fig:TorqueMap} shows the percentage change of the torque ripple depending on $J$ and $\CurrentAngle$. It becomes clear that the torque ripple is drastically reduced. Especially near to the optimized $\CurrentAngle=\SI{-11.6}{^\circ}$ and for low $J$ the torque ripple is decreased. Only outside the usual operating points $\StdTorque$ can slightly increase.

\section{Conclusion}
In this paper, we have optimized a PMSM for multiple operating points by adjusting both rotor and stator parameters and shape as well as the axial and radial dimensions by enhancing the framework from \cite{Wiesheu_2023ab}. We have shown that several objectives (cost, torque ripple, losses) can be reduced simultaneously. Despite having a large number of optimization variables, the nonlinear optimization converges quickly due to the use of gradient based optimization involving analytical expressions. For future research, including demagnetization and mechanical stresses in the optimization process are key aspects to guarantee the feasibility of the generated designs.

\section*{Acknowledgment}
The work of Michael Wiesheu and Theodor Komann is supported by the joint DFG/FWF Collaborative Research Centre CREATOR (CRC -- TRR361/F90) at TU Darmstadt, TU Graz and JKU Linz as well as the  Graduate School CE within the Centre for Computational Engineering at TU Darmstadt.

\printbibliography

\end{document}